\documentclass[12pt,leqno]{amsart}
\usepackage{amsmath}
\usepackage{amssymb}
\usepackage{amsthm} 
\usepackage{epsf}
\usepackage{graphicx}
\usepackage{esint} 
\usepackage{dsfont}
\usepackage[colorlinks,pagebackref, breaklinks]{hyperref} 
\hypersetup{backref,  pdfpagemode=FullScreen,  colorlinks=true} 

\numberwithin{equation}{section}

\theoremstyle{remark}

\newcommand{\ms}{\medskip}

\newcommand{\R}{\mathbb{R}}
\renewcommand{\H}{\mathcal H}

\newcommand{\dist}{\,\mathrm{dist}}

\newcommand{\sm}{\setminus}

\newcommand{\wt}{\widetilde}

\usepackage{color}

\definecolor{violet}{RGB}{180, 50, 180}
\definecolor{new3}{cmyk}{0, 0.61, 0.87, 0}

\begin{document}

\title{Harmonic measure on sets of codimension larger than one}
\author[David]{Guy David}
\address{Guy David. Univ Paris-Sud, Laboratoire de Math\'ematiques, UMR8628, Orsay, F-91405}
\email{guy.david@math.u-psud.fr}
\author[Feneuil]{Joseph Feneuil}
\address{Joseph Feneuil. School of Mathematics, University of Minnesota, Minneapolis, MN 55455, USA}
\email{jfeneuil@umn.edu}
\author[Mayboroda]{Svitlana Mayboroda}
\address{Svitlana Mayboroda. School of Mathematics, University of Minnesota, Minneapolis, MN 55455, USA}
\email{svitlana@math.umn.edu}
\thanks{Mayboroda is supported in part by the Alfred P. Sloan Fellowship, the NSF INSPIRE Award DMS 1344235, NSF CAREER Award DMS 1220089. David is supported in part by the ANR, programme blanc GEOMETRYA ANR-12-BS01-0014.
Feneuil is partially supported by the ANR project ``HAB'' no. ANR-12-BS01-0013.
This work was supported by a public grant as part of the Investissement d'avenir project, 
reference ANR-11-LABX-0056-LMH, LabEx LMH. 
Part of this work was completed during Mayboroda's visit to Universit\'e Paris-Sud, Laboratoire de Math\'ematiques, Orsay, and Ecole Polytechnique, PMC, and we thank the corresponding Departments and Fondation Jacques Hadamard for support and hospitality.
}
\date{}
\maketitle

\ms\noindent{\bf Abstract.}
We introduce a new notion of a harmonic measure for a $d$-dimensional set in $\R^n$ with $d<n-1$, that is, when the codimension is strictly bigger than 1. 
Our measure is associated to a degenerate elliptic PDE, it gives rise to a comprehensive elliptic theory, and, most notably, it is absolutely continuous with respect to  the $d$-dimensional Hausdorff measure 
on reasonably nice sets. This note provides general strokes of the proof of the latter statement for Lipschitz graphs with small Lipschitz constant.

\ms\noindent{\bf R\'esum\'e en Fran\c cais.}
On introduit une nouvelle notion de mesure harmonique sur un ensemble $\Gamma\subset \R^n$ 
Ahlfors-r\'egulier de dimension $d < n-1$. 
Notre mesure est associ\'ee \`a un op\'erateur diff\'erentiel lin\'eaire elliptique d\'eg\'en\'er\'e $L$, 
et a les m\^emes propri\'et\'es g\'en\'erales qu'en codimension $1$ (mesure doublante, principe
de comparaison pour les fonctions $L$-harmoniques positives). 
De plus elle est absolument continue par rapport \`a la mesure de Hausdorff de dimension $d$ 
dans des cas simples. Cette note d\'ecrit la d\'emonstration des propri\'et\'es g\'en\'erales, 
et de l'absolue continuit\'e quantifi\'ee quand $\Gamma$ est un petit graphe Lipschitzien et $L$ est bien choisi.

\ms\noindent{\bf Key words/Mots cl\'es.}
Harmonic measure in higher codimension, degenerate elliptic operators, absolute continuity, Dahlberg's theorem /
Mesure harmonique en grande codimension, op\'erateurs elliptiques d\'eg\'en\'er\'es, continuit\'e absolue, th\'eor\`eme de Dahlberg.

\ms\noindent

\tableofcontents

\section{Introduction}
\label{S1}

Recently lots of progress was made in the study of relations between regularity properties 
of the harmonic measure $\omega$ on the boundary of a domain of $\R^n$ (for instance, 
its absolute continuity with respect to the Hausdorff measure $\H^{n-1}$) 
and the regularity of the domain (for instance, rectifiability properties  of the boundary). 
In short, the emerging philosophy is that  the  
rectifiability of the boundary is necessary for the absolute continuity of $\omega$ with respect to 
$\H^{n-1}$, and that rectifiability along with suitable connectedness assumptions is sufficient. Omitting for now precise definitions, let us recall the main results in this regard. 
The celebrated 1916 theorem of F.\& M. Riesz has established the 
absolute continuity of the harmonic measure for a simply connected domain in the complex plane, with a rectifiable boundary \cite{RR}. The quantifiable analogue of this result (the
$A^\infty$ property of harmonic measure) was obtained by Lavrent'ev in 1936 \cite{Lv} 
and the local version, pertaining to subsets of a rectifiable boundary, was proved by 
Bishop and Jones in 1990 \cite{BJ}. In the latter work the authors also showed that 
some connectedness is necessary for the
absolute continuity of $\omega$ with respect to $\H^{n-1}$, for there exists a planar set with  rectifiable boundary for which the harmonic measure is singular with respect to the Lebesgue measure. The situation in higher dimensions, $n\geq 3$, is even more complicated. 
The absolute continuity of $\omega$ with respect to $\H^{n-1}$ was first established by Dahlberg 
on Lipschitz graphs \cite{Dah} and then extended to non-tangentially accessible, NTA, domains in \cite{DJ}, \cite{Se}. Roughly speaking, the non-tangential accessibility is an 
assumption of quantifiable connectedness, 
which requires the
presence of the interior and exterior corkscrew points, as well as Harnack chains. 
Similarly to the lower-dimensional case, counterexamples 
show that some topological restrictions are needed for the
absolute continuity of $\omega$ with respect to $\H^{n-1}$ \cite{Wu}, \cite{Z}. 
Much more recently, in \cite{HM1}, \cite{HMU}, \cite{AHMNT}, the authors proved that  under a (weaker) 1-sided NTA assumption, the
uniform rectifiability of the boundary is equivalent to the complete set of NTA  conditions 
and hence, is equivalent to the
absolute continuity of harmonic measure with respect to the Lebesgue measure. 
Finally, in 2015 the full converse, ``free boundary"  
result was obtained and established that rectifiability is necessary for  the
absolute continuity of harmonic measure with respect to $\H^{n-1}$ in any dimension $n\geq 2$ (without any additional topological assumptions) \cite{AHM3TV} .  

The purpose of this work is to start the investigation
of similar properties for domains with a lower-dimensional boundary $\Gamma$.  

Here we shall systematically assume that $\Gamma$ is Ahlfors-regular of some dimension
$d < n-1$, which for the moment does not need to be an integer. 
This means that there is a constant $C_0 \geq 1$ such that
\begin{equation} \label{1.1}
C_0^{-1} r^d \leq \H^d(\Gamma \cap B(x,r)) \leq C_0 r^d
\ \text{ for } x\in \Gamma \text{ and } r > 0.
\end{equation}
We want to define an analogue of harmonic measure, that will be defined on $\Gamma$ and
associated to a divergence form operator on $\Omega = \R^n \sm \Gamma$. We still write
the operator as $L = \rm{div} A \nabla$, with $A : \Omega \to \mathbb{M}_n(\R)$, 
and we write the ellipticity condition with a different homogeneity, i.e., require that for some 
$C_1 \geq 1$, 
\begin{eqnarray} \label{1.2.1}
&& \dist(x,\Gamma)^{n-d-1} A(x)\xi \cdot \zeta \leq C_1  |\xi| \,|\zeta|
\ \text{ for } x\in \Omega \text{ and } \xi, \zeta \in \R^n, \\[4pt]
\label{1.2.2}
&& \dist(x,\Gamma)^{n-d-1} A(x)\xi \cdot \xi \geq C_1^{-1}  |\xi|^2
\ \text{ for } x\in \Omega \text{ and } \xi \in \R^n.
\end{eqnarray}
The effect of this normalization should be to incite the analogue of the Brownian motion here
to get closer to the boundary with the right probability;  for instance if $\Gamma = \R^d \subset \R^n$
and $A(x) = \dist(x,\Gamma)^{-n+d+1} I$, it turns out that the effect of $L$ on functions $f(x,t)$ 
that are radial in the second variable $t\in \R^{n-d}$ is the same as for the Laplacian on $\R^{d+1}_+$. 
In some sense, we create Brownian travelers  which treat $\Gamma$ as a ``black hole": they detect more mass and they are more attracted to $\Gamma$ than a standard Brownian traveler governed by the Laplacian would be.

With merely these assumptions, we solve a first Dirichlet problem for $Lu = 0$, 
prove the maximum principle, the De Giorgi-Nash-Moser estimates and the 
Harnack inequality for solutions, use this to define a harmonic measure associated to $L$, show that it is doubling, 
and prove the comparison principle for positive $L$-harmonic functions that vanish at the boundary 
(see below).

Then we take stronger assumptions, both on the geometry of $\Gamma$ and the choice of $L$,
and try to prove that the harmonic measure is absolutely continuous with respect to 
$\H^d_{\vert \Gamma}$. Here we assume that $d$ is an integer and $\Gamma$ is the graph 
of a Lipschitz function $F : \R^d \to \R^{n-d}$, with a small enough Lipschitz constant.
As for $A$, we assume that $A(x) = D(x)^{-n+d+1} I$ for $x\in \Omega$, with 
\begin{equation} \label{1.3}
D(x) = \Big\{ \int_\Gamma |x-y|^{-d-\alpha} d\H^d(y) \Big\}^{-1/\alpha}
\end{equation}
for some constant $\alpha > 0$. Notice that because of \eqref{1.1}, $D(x)$ is equivalent to 
$\dist(x,\Gamma)$; when $d=1$ we can also take $A(x) = \dist(x,\Gamma)^{-n+d+1} I$, but
when $d \geq 2$ we do not know whether $\dist(x,\Gamma)$ is smooth enough to work. 
In \eqref{1.3}, we could also replace $\H^d$ with another Ahlfors-regular measure on $\Gamma$. 

With these additional assumptions we prove that the harmonic measure described above 
is absolutely continuous with respect to $\H^d_{\vert \Gamma}$, with a density which is 
a Muckenhoupt $A_\infty$ weight. In other words, we establish an analogue of Dahlberg's result \cite{Dah} for domains with a higher co-dimensional boundary given by a Lipschitz graph 
with a small Lipschitz constant.
It is not so clear what is the right condition for this in terms of $A$, 
but the authors still hope that 
a good condition on $\Gamma$ is its uniform rectifiability. 
Notice that in remarkable contrast with the case of codimension 1, we do not state an additional quantitative connectedness condition on $\Omega$, 
such as the Harnack chain condition in codimension $1$; this is because 
such conditions are automatically satisfied when 
$\Gamma$ is Ahlfors-regular with a large codimension. 
See Harnack's inequality below.

\section{Sketch of proofs of the main results}

The authors intend to give detailed proofs of the results above somewhere else \cite{DFM},
but think that in the mean time the rapid description of the proof given below will give a fair idea of the arguments.
We first introduce some notation. Set $\delta(x) = \dist(x,\Gamma)$
and $w(x) = \delta(x)^{-n+d+1}$ for $x\in \Omega = \R^n \sm \Gamma$,
and denote by $\sigma$ the restriction to $\Gamma$ of $\H^d$. 
Denote by $W=\dot W^{1,2}_w(\Omega)$ the weighted Sobolev space of functions 
$u \in L^1_{loc}(\Omega)$
whose distribution gradient in $\Omega$ lies in $L^2(\Omega,w)$, and set
$||u||_W = \big\{\int_\Omega |\nabla u(x)|^2 w(x) dx \big\}^{1/2}$ for $f\in W$.
Finally denote by $H$ or $\dot H^{1/2}(\Gamma)$ the set of measurable functions $g$ on $\Gamma$
for which $||g||_H^2 = \int_\Gamma\int_\Gamma {|g(x)-g(y)|^2 \over |x-y|^{d+1}}d\sigma(x) d\sigma(y)$ is finite.

Before we solve Dirichlet problems we construct two bounded linear operators
$T : W \to H$ (a trace operator) and $E : H \to W$ (an extension operator), such that
$T \circ E = I_H$. The trace of $u \in W$ is such that
$Tu(x) = \lim_{r \to 0} \fint_{B(x,r)} u(y) dy :=  \lim_{r \to 0} {1 \over |B(x,r)|} \int u(y) dy$,
and even, analogously to the Lebesgue density property,
$\lim_{r \to 0} \fint_{B(x,r)} |u(y)-Tu(x)| dy = 0$ for $\sigma$-almost every $x\in \Gamma$.
Similarly, we check that if $g\in H$, then 
$ \lim_{r \to 0} \fint_{\Gamma \cap B(x,r)} |g(y)-g(x)| d\sigma(y) = 0$
for $\sigma$-almost every $x\in \Gamma$. The proofs are standard and easy;
we typically use the fact that $|u(x)-u(y)| \leq \int_{[x,y]} |\nabla u|$ for almost all choices
of $x$ and $y \in \Omega$, for which we can use the absolute continuity of $u\in W$ on 
(almost all) line segments, plus the important fact that, by \eqref{1.1}, 
$\Gamma \cap \ell = \emptyset$ for almost every line $\ell$
(see the discussion about Harnack below for an even better result). For the construction
of $E$, we start with the standard proof of the Whitney extension theorem.

Once we have these operators, we easily deduce from the Lax-Milgram theorem that for $g\in H$, 
there is a unique weak solution $u\in W$ of $Lu=0$ such that $Tu = g$. 
For us a weak solution is a function $u\in W$ such that
$\int_\Omega A(x)\nabla u(x) \cdot \nabla \varphi(x) dx = 0$ for 
$\varphi \in C_0^\infty (\Omega)$, the space of infinitely differentiable functions which are compactly 
supported in $\Omega$.

Then we follow the Moser iteration scheme to study the weak solutions of $Lu=0$, as we would do
in the standard elliptic case in codimension $1$.  
This leads to the quantitative boundedness (a.k.a. Moser bounds) and the quantitative H\"older continuity
(a.k.a. De Giorgi-Nash estimates), in  an interior or boundary ball $B$,
of any weak solution of $Lu = 0$ in $2B$ such that $Tu = 0$ on $\Gamma \cap 2B$ 
when the intersection is non-empty.  The boundary estimates are trickier, because we do not have the conventional ``fatness" of the complement of the domain, and it is useful to know that  
\begin{equation} \label{1.4}
\fint_{B(x,r)} |u(y)| dy \leq C r^{-d} \int_{B(x,r)} |\nabla u(y)| w(y) dy
\end{equation} 
for $u \in W$, $x\in \Gamma$, and $r>0$ such that $Tu = 0$ on $\Gamma \cap B(x,r)$
and that, if $V(x,r)$ denotes $\int_{B(x,r)} w(y) dy$,
\begin{equation} \label{1.5}
\Big\{\frac1{V(x,r)}\int_{B(x,r)} \Big|u(y) - \fint_{B(x,r)} u \Big|^p w(y) dy\Big\}^{1/p}
\leq C r \Big\{ \frac{1}{V(x,r)}\int_{B(x,r)} |\nabla u(y)|^2 w(y) dy \Big\}^{1/2}
\end{equation}
for $u\in W$,  $x\in \overline\Omega=\R^n$, $r>0$, and 
$p\in \left[1,\frac{2n}{n-1}+\delta'\right)$ for some $\delta'>0$. 
For \eqref{1.4} we write that $|u(y)-u(\xi)| \leq \int_{[y,z]} |\nabla u|$ for almost all $y\in B(x,r)$, 
$z\in \Gamma \cap B(x,r)$ and almost all $\xi \in B(z,\varepsilon r)$, and then integrate 
and take a limit; for the weighted Poincar\'e inequality \eqref{1.5} we observe that $\Gamma$ 
does not interfere with the proof of Poincar\'e inequality with weights in \cite{FKS}, 
because almost every line of $\R^n$ is contained in $\Omega$.

We also need to find Harnack chains, and in fact there exists a constant $c>0$, that depends only
on $C_0$, $n$, and $d < n-1$, such that for $\Lambda \geq 1$ and $x_1, x_2 \in \Omega$ 
such that $\dist(x_i,\Gamma) \geq r$ and $|x_1-x_2| \leq \Lambda r$, we can find two points
$y_i \in B(x_i,r/2)$ such that $\dist([y_1,y_2],\Gamma) \geq c \Lambda^{-d/(n-d-1)} r$.
That is, there is a thick tube in $\Omega$ that connects the two $B(x_i,r/2)$.

With all these ingredients, we can follow the standard proofs for elliptic divergence form operators, prove a suitable version of the Harnack inequality and the maximum principle, 
solve the Dirichlet problem for continuous functions with compact support on $\Gamma$, 
define harmonic measures $\omega^x$ for $x\in \Omega$ (so that $\int_\Gamma g d\omega^x$ 
is the value at $x$ of the solution of the Dirichlet problem for $g$), prove that $\omega^x$
is doubling, and even establish the comparison principle that says that if $u$ and $v$ are
positive weak solutions of $Lu= 0$ such that $Tu = Tv = 0$ on $B(x,2r) \cap\Gamma$, 
with $x\in \Gamma$, then $u$ and $v$ are comparable in $B = B(x,r)$, i.e.,
$\sup_{z\in B \sm \Gamma} (u(z)/v(z)) \leq C \inf_{z\in B \sm \Gamma} (u(z)/v(z))$.

Let us point out that while the invention of a harmonic measure which serves the higher co-dimensional boundaries, which is associated to a linear PDE, and which is absolutely continuous with respect to the Lebesgue measure on reasonably nice sets, 
is the main focal point of this work, various versions of degenerate elliptic operators and weighted Sobolev spaces have of course appeared in the literature over the years. 
Some versions of some of the results 
listed above or similar ones can 
be found, e.g., in \cite{Ancona}, \cite{FKS}, \cite{Hajlasz}, \cite{NDbook}, and other sources, not to mention that modulo sorting out some trace questions
we followed the traditional line of development. Since we did not rely on previous work, we hope to be forgiven for not reviewing the corresponding literature. In our setting, a more intricate and difficult part of the story is certainly the proof of the $A^\infty$ property that we discuss below.

\ms
We turn to the proof of absolute continuity (with $A_\infty$ estimates) in the case of small
Lipschitz graphs, and for $L = \rm{div} A \nabla$, with $A(x) = D(x)^{d+1-n} I$. 
The main initial ingredient is a bilipschitz change of variable $\rho : \R^n \to \R^n$, such that 
$\rho(\R^d) = \Gamma$, and which is also often nearly isometric in the last $n-d$ variables. 
That is, write $z = (x,t)$ the running point of $\R^n$, with $x\in \R^d$ and $t\in \R^{n-d}$;
we construct $\rho$, essentially by hand, so that in particular 
$\displaystyle\left|{D(\rho(x,t)) \over |t|} - 1 \right|^2 {dx dt \over |t|^{n-d}}$ is a Carleson measure
relative to $\R^d \subset \R^n$.  

The interest of this $\rho$ is that, when we conjugate $L$ with $\rho$, we get an
operator $\wt L = \rm{div} \wt A \nabla$ on $\wt \Omega = \R^n \sm \R^d$, and 
with our good choices of distance $D$, we can write $\wt A$ as a sum
$\wt A_1 + \wt A_2$, where $\wt A$ and $\wt A_1$ satisfy the same ellipticity conditions \eqref{1.2.1} and \eqref{1.2.2}
as $A$ (but now relative to $\R^d$), and in addition each $\wt A_1(x,t)$ is a block matrix
$\left(\begin{array}{cc} A_1'&0 \\ 0 & A_1'' \end{array}\right)$,
with $A_1'\in \mathbb{M}_{d}(\R)$, while $A_1'' \in \mathbb{M}_{n-d}(\R)$ is a multiple of the identity, $A_1''(x, t)=a_1''(x, t) I$, and the scalar function $a_1''(x,t)$ is such that
$\big||t|\nabla a_1''(x,t)\big|^2 {dx dt \over |t|^{n-d}}$ is a Carleson measure. 
As for $\wt A_2$, it is small in the sense that 
$\big| \wt A_2(x,t) \big|^2 {dx dt \over |t|^{n-d}}$ 
is a Carleson measure. Of course the desired $A_\infty$ result for $L$ and $\Gamma$ follows from
the same thing for $\wt L$ and $\R^d$.

We prove the $A_\infty$ bounds in two steps. 
The first one uses integrations by parts to prove that whenever the matrix of coefficients $\wt A$ satisfies conditions for the structure and smoothness/size described above, 
the following localized bound for the square function
of an $L^2$-solution 
of $\wt L u = 0$, in terms of the maximal nontangential function holds: 
for every cube $Q\subset \R^d$
\begin{equation} \label{1.7}
\|S^Q u\|_{L^2(Q)}^{2} \lesssim \|\mathcal{N}^Q u\|_{L^2(2Q)}^{2} + \int_Q u^2\, dx.
\end{equation}
Here, the localized square function is
$$S^Q u(x, 0):=\left(\int_{(y,s)\in \gamma^Q (x,0)} |\nabla u(y,s)|^2 \, \frac{dyds}{|(y,s)-(x,0)|^{n-2}}\right)^{1/2}, \quad x\in \R^d, $$
the localized non-tangential maximal function is
$$
{\mathcal N}^Q u(x, 0):=\sup_{(y,s)\in \widetilde \gamma^Q (x,0)} |\nabla u(y,s)|, \quad x\in \R^d, $$
and $\gamma^Q$ is a truncated cone 
$$\gamma^Q(x,0)=\{(y,s)\in \R^n\setminus \R^d:\, |y-x|\leq |s|, \,0<|s|< l(Q)\}, $$
while $\widetilde \gamma^Q$ is a truncated cone with bigger aperture, that is, the condition $|y-x|\leq |s|$ in the definition of $\gamma^Q$ is substituted by  $|y-s|<C |s|$ for $C$ large enough but depending on $n, d$ only.
By the maximum principle \eqref{1.7} entails that if $u$ is a solution $\wt L u = 0$, whose boundary values on $\Gamma$ are the characteristic function of a Borel set, then 
$\left||t| \nabla u(x,t)\right|^2 {dx dt \over |t|^{n-d}}$ is a Carleson measure.
The proof of \eqref{1.7} builds on some ideas of \cite{KePiDrift}, though the fact that there are $n-d$ linearly independent vectors in the transversal direction to the boundary brings some completely new difficulties.

We complete the proof with an argument resonating with \cite{KKPT}, \cite{KKiPT} 
which shows the following. For any elliptic operator $L$ satisfying \eqref{1.2.1}--\eqref{1.2.2} 
and  $\Gamma=\R^d$, if all the $L$-harmonic extensions of characteristic functions of 
Borel subsets of $\R^d$ satisfy the property that 
$\left||t| \nabla u(x,t)\right|^2 {dx dt \over |t|^{n-d}}$ is a Carleson measure, 
then  the harmonic measure for $L$ is $A_\infty$-equivalent the Lebesgue measure on $\R^d$. 
Note that at this stage no structural or smoothness conditions on the underlying matrix are required. 
It is the proof of \eqref{1.7} and the core Carleson measure estimates on solutions that relies on the special structure of the operator. Recall though that even in codimension 1 one does not expect 
the absolute continuity of the harmonic measure with respect to the Lebesgue measure to hold for all elliptic matrices, due to the counterexamples in \cite{CFK}.

\end{document}